\documentclass[10pt]{amsart}

\usepackage{amssymb,amsmath,amsthm}
\usepackage{graphicx}

\theoremstyle{plain}

\usepackage{graphicx,tikz}
\newtheorem{theorem}{Theorem}

\theoremstyle{definition}

\theoremstyle{remark}

\begin{document}

\title[]{ Electrostatic Interpretation of Zeros\\of orthogonal polynomials   }

\author[]{Stefan Steinerberger}
\address{Department of Mathematics, Yale University}
\email{stefan.steinerberger@yale.edu}

\begin{abstract}
We study the differential equation
$ - (p(x) y')' + q(x) y' = \lambda y,$
where $p(x)$ is a polynomial of degree at most 2 and $q(x)$ is a polynomial of degree at most 1. This includes the classical Jacobi polynomials, Hermite polynomials, Legendre polynomials, Chebychev polynomials and Laguerre polynomials. We
provide a general electrostatic interpretation of zeros of such polynomials: a set of distinct, real numbers $\left\{x_1, \dots, x_n\right\}$ satisfies
$$    p(x_i) \sum_{k = 1 \atop k \neq i}^{n}{\frac{2}{x_k - x_i}} = q(x_i) -  p'(x_i) \qquad \mbox{for all}~ 1\leq i \leq n$$
if and only if they are zeros of a polynomial solving the differential equation. We also derive a system of ODEs depending on $p(x),q(x)$ whose solutions converge to the zeros of the orthogonal polynomial at an exponential rate.
\end{abstract}

\maketitle

\vspace{-20pt}

\section{Introduction}
\subsection{Introduction.}
We start by describing an 1885 result of Stieltjes for Jacobi polynomials \cite{stieltjes}. Jacobi polynomials $P_n^{\alpha, \beta}(x)$, for real $\alpha, \beta > -1$, are the unique (up to a constant factor) solutions of the equation
$$ (1 - x^2 ) y''(x) - \left( \beta - \alpha - (\alpha + \beta + 2)x\right)y'(x)  = n(n + \alpha + \beta + 1) y(x)$$
where $n \in \mathbb{N}$. The solution is a polynomial of degree $n$ having all its zeros in $(-1,1)$.  Stieltjes defined a notion of energy for any set $\left\{x_1, \dots, x_n \right\} \subset (-1,1)$ as 
$$ E = - \sum_{i,j = 1 \atop i \neq j}^{n}{\log{|x_i - x_j|}} - \sum_{i=1}^{n}{\left( \frac{\alpha+1}{2}\log{|x-1|} + \frac{\beta+1}{2}\log{|x+1|}\right)}$$
and showed that the minimal energy is exactly given by the zeros of the Jacobi polynomial. Differentiating this expression in all the $n$ variables leads to
a system of equations describing an electrostatic equilibrium
$$ \sum_{k=1 \atop k \neq i}^{n}{ \frac{1}{x_k - x_i} } =  \frac{1}{2}\frac{\alpha+1}{x_i - 1} + \frac{1}{2} \frac{\beta + 1}{x_i + 1} \qquad \mbox{for all}~1 \leq i \leq n.$$
A traditional application of the formula is to establish monotonicity of zeros with respect to the parameters $\alpha, \beta$. This argument can be found in Szeg\H{o}'s book \cite{szeg} who also discusses Laguerre and Hermite polynomials. These types of questions have been studied by a large number of people, we refer to \cite{ahmed, dimitrov, hendriksen, ismail}, the survey \cite{survey},
the 1978 survey and 2001 book by F. Calogero \cite{calogero, calogero2}, the more recent papers of F. A. Gr\"unbaum \cite{grunbaum, grunbaum2} 
and M. E. H. Ismail \cite{ismail, ismail2} and references therein.

\subsection{Equilibrium.}
The purpose of this paper is to prove a simple general result that characterizes zeros of orthogonal polynomials. The statement
itself is a bit more general but most interesting when applied to the classical orthogonal polynomials. We first discuss
 second order equations, then describe system of ODEs converging to zeros and conclude with a short remark on higher order equations.

\begin{theorem}  Let $p(x), q(x)$ be polynomials of degree at most 2 and 1, respectively. Then the set $\left\{x_1, \dots, x_n\right\}$, assumed to be in the domain of definition, satisfies
$$   p(x_i)  \sum_{k = 1 \atop k \neq i}^{n}{\frac{2}{x_k - x_i}} = q(x_i) -  p'(x_i) \qquad \mbox{for all} \quad 1\leq i \leq n$$
if and only if 
$$y(x) = \prod_{k=1}^{n}{(x-x_k)} \quad \mbox{solves} \quad - (p(x) y')' + q(x) y' = \lambda y \quad \mbox{for some} ~~ \lambda \in \mathbb{R}.$$
\end{theorem}
In classical applications, there is a unique polynomial solution of fixed degree corresponding to a fixed value of $\lambda$. This removes $\lambda$ as a variable and leads to a complete characterization.
One direction of the statement (zeros of a polynomial solutions satisfy the system of equations) is a fairly straightforward computation and the outline of the argument can be found, for example, in \cite{grunbaum} or, as a remark, in \cite[\S 3]{ismail}. Moreover, results in this direction can be obtained for much more general differential equations of higher order (see, for example,  \cite{calogero, calogero2,avila}). Arguments in the other direction seem to exist only in special cases, as in the work of Stieltjes, and are based on interpreting the system of equations as the critical point of an associated energy functional (see e.g. \cite{ismail, survey}). Our argument proceeds by a different route and bypasses considerations of an underlying energy.\\

We remark that there is no assumption of orthogonality nor is there any restriction on the domain which may be either bounded or unbounded but we do assume that the $n$ points are contained in the domain of definition. Let us quickly consider two special cases: we start with the Hermite differential equation
$$ -y'' + x y' = \lambda y.$$
This corresponds to $p(x) \equiv 1$ and $q(x) \equiv x$. We deduce from Theorem 1 that the zeros of the $n-$th solution satisfy a relationship also going back to Stieltjes, see \cite{ahmed1},
$$     \sum_{k = 1 \atop k \neq i}^{n}{\frac{2}{x_k - x_i}} = x_i \qquad \mbox{for all}~1 \leq i \leq n.$$
Returning to a special case of Jacobi polynomials $\alpha = 0 = \beta$ (for simplicity of exposition), we obtain $p(x) = x^2 -1$, $q(x) = 0$ and thus 
$$   (x_i^2-1)  \sum_{k = 1 \atop k \neq i}^{n}{\frac{2}{x_k - x_i}} = 2x_i \qquad \mbox{for all} \quad 1\leq i \leq n$$
which is easily seen to be equivalent to Stieltjes' electrostatic equilibrium.

\subsection{A System of ODEs.} One interesting byproduct of our approach is a procedure that allows for the computation of zeros without ever computing the polynomial.
More precisely, for any given set $x_1(0) < x_2(0) < \dots < x_n(0)$ (assumed to be in the domain of definition), we consider the system of ordinary differential equations
$$    \frac{d}{dt} x_i(t) =  -p(x_i) \sum_{k = 1 \atop k \neq i}^{n}{\frac{2}{x_k(t) - x_i(t)}} +  p'(x_i(t)) -  q(x_i(t))   \qquad \qquad (\diamond)$$
The result above implies that the unique stationary point of this system of ODEs is given by the $n$ zeros of a polynomial solution of degree $n$. 
We will now show that under standard assumptions on the equation $ - (p(x) y')' + q(x) y' = \lambda y,$ the system converges to this fixed point at an exponential rate.
We require somewhat stronger assumptions than we do for Theorem 1 and ask that
\begin{enumerate}
\item the equation  $ - (p(x) y'(x))' + q(x) y'(x) = \lambda y(x)$ is a Sturm-Liouville problem with a discrete set of solutions indexed by $0 = \lambda_0 < \lambda_1 < \lambda_2 < \dots$  and
\item the solution corresponding to $y_n$ is a polynomial of degree $n$ for all $n \in \mathbb{N}$
\end{enumerate}
These assumptions cover the classical polynomials but are fairly strong, we refer to Bochner's theorem and its nice exposition in \cite[\S 20]{ismail3}.
\begin{theorem}
The system $(\diamond)$ converges for all initial values $x_1(0) < \dots < x_n(0)$ to the zeros $x_1 < \dots < x_n$ of the degree $n$ polynomial solving the equation. Moreover,
$$ \max_{1 \leq i \leq n}{ |x_i(t) - x_i|} \leq c e^{- \sigma_n t},$$
where $c > 0$ depends on everything and $\sigma_n \geq \lambda_n - \lambda_{n-1}$.
\end{theorem}
We are not aware of any result of this type being known; the connection between partial differential equations and zeros of polynomials (or poles of rational functions) is, of course classical, we refer to the extensive
survey of Calogero \cite{calogero} and his more recent book \cite{calogero2}. Our underlying dynamical system is the \textit{backward} heat equation (well-posed for algebraic reasons) and existence
of a solution for all time is a consequence of Sturm-Liouville theory.
We assume that the initial values are all contained in the domain where the solution is defined. For classical Sturm-Liouville problems, the Weyl law suggests 
$\lambda_{n} - \lambda_{n-1} \sim n$. The rapid convergence can be easily observed in examples.
The assumptions are satisfied for all classical polynomials. We consider, as a specific example, the Laguerre polynomials satisfying
$$ - x y'' + (x-1)y' = ny $$
or, in our notation, $p(x) = x = q(x)$. Theorem 2 then implies that for any choice of initial values $x(0) < y(0) < z(0)$, the system of ordinary differential equations 
\begin{align*}
\dot x(t) &= \frac{2x(t)}{x(t) - y(t)} + \frac{2x(t)}{x(t) - z(t)} +1 - x(t) \\
\dot y(t) &= \frac{2y(t)}{y(t) - x(t)} + \frac{2y(t)}{y(t) - z(t)} +1 - y(t) \\
\dot z(t) &= \frac{2z(t)}{z(t) - x(t)} + \frac{2z(t)}{z(t) - y(t)} +1 - z(t) 
\end{align*}
has a solution for all $t > 0$.  As $t \rightarrow \infty$, the solutions converge to the zeros of the third Laguerre polynomial
$$ L_3(x) = x^3 - 9x^2 + 18x - 6$$
and thus
$$ x(t) \rightarrow 0.4157\dots,~ y(t) \rightarrow 2.2942\dots \quad \mbox{and} \quad z(t) \rightarrow 6.2899\dots.$$
Moreover, this convergence happens at an exponential rate (roughly with speed $e^{-t}$).
We consider a large-scale example (see Fig. 1) next. The equation
$$- \frac{d}{dx}\left( (1-x^2) \frac{d}{dx}y(x)\right) = n(n+1) y(x)$$
is solved by the Legendre polynomials $P_n$ defined on $(-1,1)$. In our notation, we have $p(x) = 1- x^2$ and $q(x) = 0$. We simulate
the system of ODEs
 $$    \frac{d}{dt} x_i(t) =  -(1-x_i(t))^2 \sum_{k = 1 \atop k \neq i}^{n}{\frac{2}{x_k(t) - x_i(t)}} - 2x_i(t)$$
for $n=100$ equations. The initial values $\left\{x_1(0), \dots, x_{100}(0)\right\}$ are chosen as uniform random variables
in the interval $(-0.1, 0.1)$. Figure 1 shows the evolution of the system and convergence to the zeros of the
Legendre polynomial of degree 100. 
\vspace{-10pt}
\begin{center}
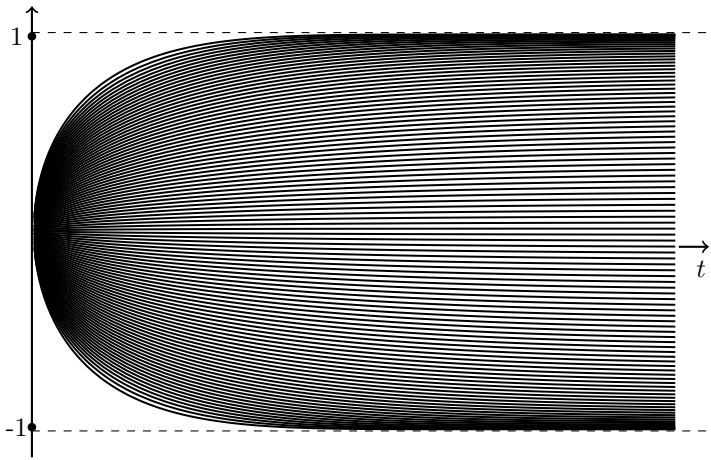
\begin{figure}[h!]
\begin{tikzpicture}[scale=1]
\draw [thick, ->] (-4.3,-3) -- (-4.3,3);
\draw [thick, ->] (4.3,-0.2) -- (4.7,-0.2);
\node at (4.6, -0.5) {$t$};
\filldraw (-4.3, 2.6) circle (0.05cm);
\draw [dashed] (-4.3,2.65) -- (4.7,2.65);
\draw [dashed] (-4.3,-2.65) -- (4.7,-2.65);
\node at (-4.5, 2.6) {1};
\filldraw (-4.3, -2.6) circle (0.05cm);
\node at (-4.5, -2.6) {-1};
\node at (-.010,0) {\includegraphics[width=0.7\textwidth]{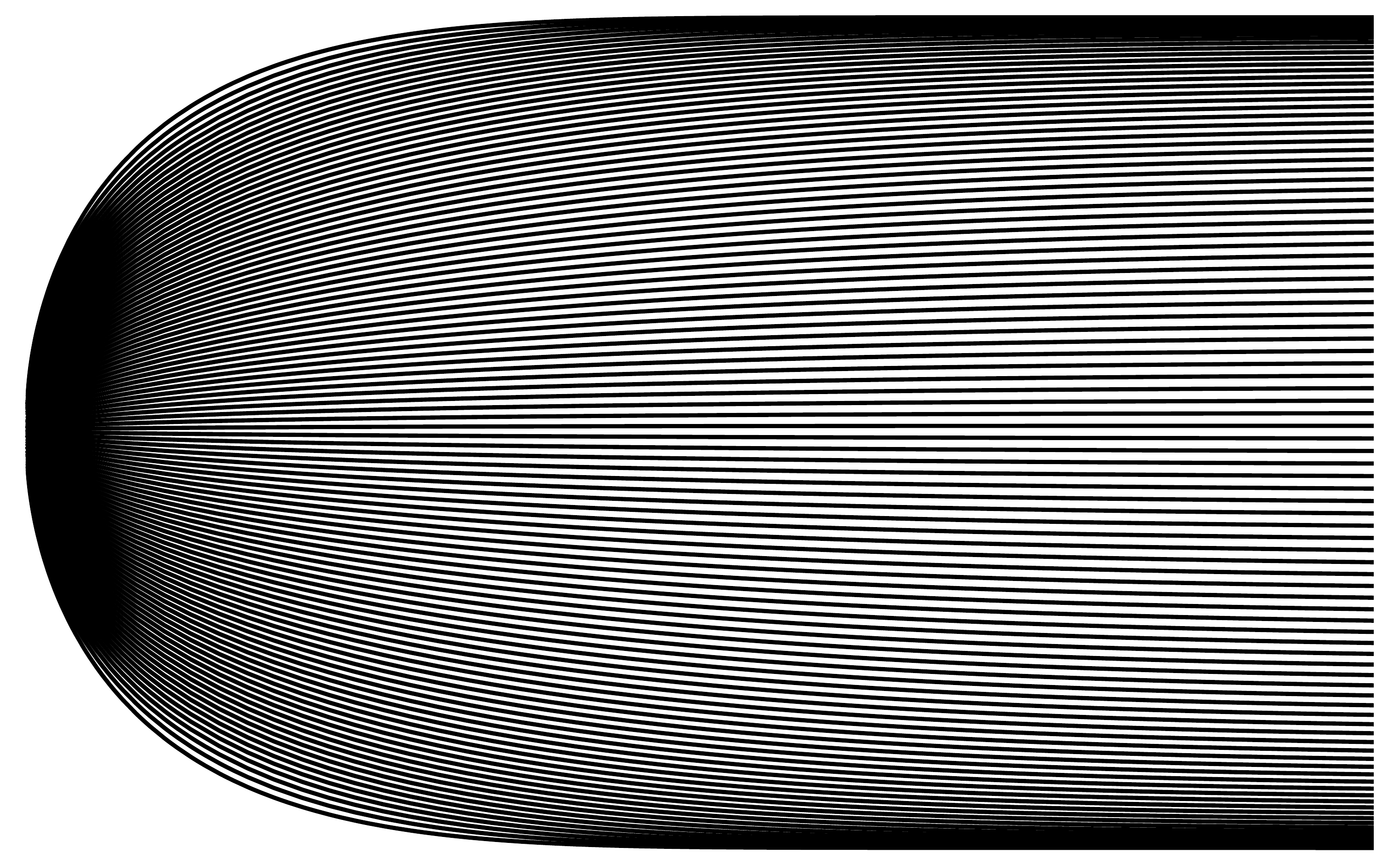}};
\end{tikzpicture}
\caption{Evolution of the system of ODEs for $0 \leq t \leq 0.01$ approaches the zeros of the Legendre polynomial $P_{100}$ in $(-1,1)$.}
\end{figure}
\end{center}
\vspace{-20pt}

As indicated in the Theorem, the constant determining the speed of exponential
convergence grows linearly in $n$: this is reflected in the rather short time-scale $(0 \leq t \leq 0.01)$ in the picture.

\subsection{Higher order equations} Theorem 1 as well as its proof immediately generalizes to polynomial solutions of the equation
$$ a_n(x) y^{(n)}(x) + a_{n-1}(x) y^{(n-1)}(x) + a_{1}(x)y'(x)  = \lambda y(x),$$
where $a_n$ is a polynomial of degree at most $n$ and the solution $y$ is assumed to only have single zeros: one half 
of the statement can be found in \cite[Proposition 1]{avila}, the other direction follows from our approach.
We are not aware of an extension of Theorem 2 since our
proof relies on Sturm-Liouville theory.

\section{Proofs}
\subsection{Proof of Theorem 1.}
\begin{proof}
One direction of the argument, showing that the zeros of any polynomial solution of the equation is in the desired electrostatic equilibrium, is a simple computation and known in greater generality \cite{ahmed1,ahmed}. However, it also follows immediately from our argument and is, for completeness sake, sketched at the end of the proof.
We assume that the system of equations is satisfied and will prove that the associated polynomial solves the differential equation. Fix $x_1, \dots, x_n$
and consider the candidate polynomial
$$ f(x) = \prod_{k=1}^{n}{(x - x_k)}.$$
We want to show that this polynomial satisfies the differential equation
$$- (p(x) f')' + q(x) f' = \lambda f \qquad \mbox{for some}~\lambda \in \mathbb{R}.$$
We introduce the function
\begin{align*}
\frac{\partial}{\partial t} u(t,x) &=  \frac{\partial}{\partial x}\left( p(x) \frac{\partial}{\partial x} u(t,x) \right) - q(x) \frac{\partial}{\partial x} u(t,x) \\
u(0,x) &= f(x)
\end{align*}
This is a parabolic partial differential equation. Generally, unless we specify the sign of $p$, it might be an inverse heat equation and the solution of such an equation need not even exist for a short amount of time. Here, however, since $p$ is at most a polynomial
of degree 2 and $q$ is a polynomial of degree at most 1, we see that the right-hand side is always a polynomial of degree at most $n$. In particular, we can rewrite the partial 
differential equation as a linear system of $n+1$ ordinary differential equations and this guarantees existence for all $t>0$.
Suppose $f(x)$ is not a solution of the differential equation. Then 
$$u_t(0,x) =  \frac{\partial}{\partial x}\left( p(x) \frac{\partial}{\partial x} f(x) \right) - q(x) \frac{\partial}{\partial x} f(x)$$
 is a polynomial of degree at most $n$ and not identically 0 (since otherwise $f$ would solve the differential equation for $\lambda = 0$). We observe that if this polynomial vanishes in
exactly $\left\{x_1, \dots, x_n\right\}$, then it has to be a multiple of $f$ and we have obtained the desired result for some $\lambda \neq 0$. If this is not the case, then $u_t(x_i,0) \neq 0$ for at least one $1 \leq i \leq n$.
We fix this value of $i$. Moreover, we note 
$$f'(x_i) =  \prod_{k=1 \atop k \neq i}^{n}{(x_i - x_k)} \neq 0.$$
The implicit function theorem implies that there is a neighborhood of 0 for which there is a function $x_i(t)$ such that
$$ u( t, x_i(t)) = 0 \qquad \mbox{and} \qquad \frac{\partial}{\partial t} x_i(t) \big|_{t = 0} \neq 0.$$ 
We now compute the expression. Differentiation implies
$$ 0 = \frac{\partial}{\partial t} u(t,x_i(t))\big|_{t=0} = u_x(0,x_i)\left(\frac{\partial}{\partial t} x_i(t) \big|_{t=0}\right)  + u_t(0,x_i).$$
We are interested in the first term, already computed the second term $u_x(x_i,0) = f'(x_i)$ and thus only need to compute the third term. A simple computation shows
\begin{align*}
u_t(t,x_i)\big|_{t=0} &=  \frac{\partial}{\partial x}\left( p(x) \frac{\partial}{\partial x} f(x) \right) - q(x) \frac{\partial}{\partial x}  f(x) \big|_{x=x_i} \\
&=  p(x) \frac{\partial^2}{\partial x^2} f(x) + (p'(x) - q(x))\frac{\partial}{\partial x}  f(x) \big|_{x=x_i}.
\end{align*}
The first term simplifies to
$$ \frac{\partial^2}{\partial x^2} f(x) \big|_{x=x_i} = 2 \sum_{k =1 \atop k \neq i}^{n}{ \prod_{j = 1 \atop  j \notin \left\{i, k\right\}}^{n}{(x_i - x_j)}}$$
and altogether we obtain
$$ 0 \neq \frac{\partial}{\partial t} x_i(t) \big|_{t = 0} = p(x_i)  \sum_{k = 1 \atop k \neq i}^{n}{\frac{2}{x_k - x_i}} + p'(x_i) - q(x_i)$$
which is a contradiction. Conversely, if $f$ is indeed a solution of the equation, then 
$$ u(t,x) = e^{\lambda t} f(x) \qquad \mbox{and thus} \qquad \frac{\partial}{\partial t} x_i(t) \big|_{t = 0} = 0$$
for all $1 \leq i \leq n$ which implies that the equations are satisfied.
\end{proof}

\begin{center}
\begin{figure}[h!]
\begin{tikzpicture}[scale=1]
\draw [thick, ->] (-4.3,-3) -- (-4.3,3);
\node at (4.7, -2.8) {$t$};
\filldraw (-4.3, 2.6) circle (0.05cm);
\draw [thick, ->] (-4.3,-2.6) -- (4.7,-2.6);
\node at (-4.7, 2.6) {371};
\filldraw (-4.3, -2.6) circle (0.05cm);
\filldraw (4.1, -2.6) circle (0.05cm);
\node at (4.1, -2.8) {4};
\node at (-4.1, 2.9) {$x$};
\node at (-4.6, -2.6) {0};
\node at (-0.02,0.04) {\includegraphics[width=0.7\textwidth]{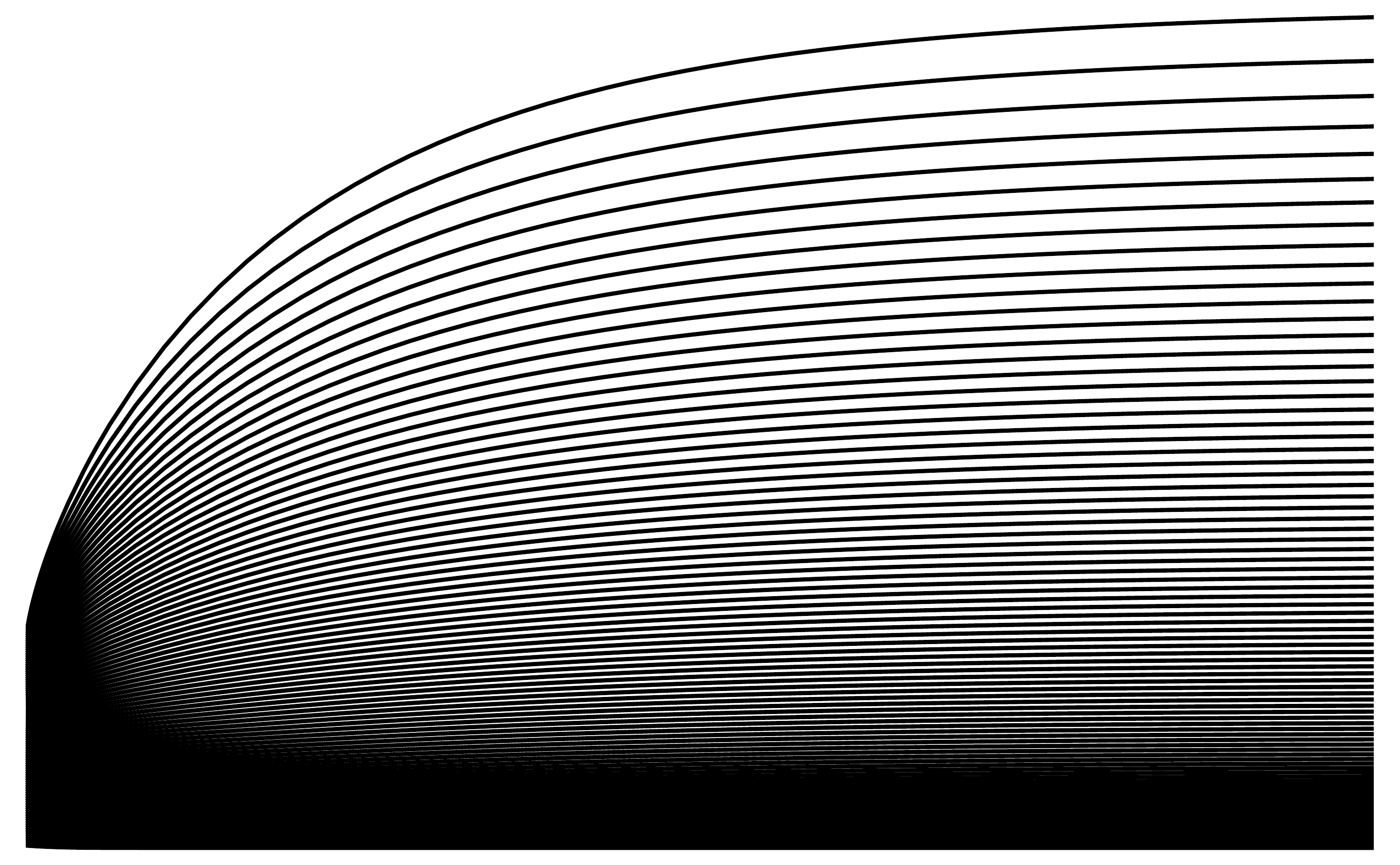}};
\end{tikzpicture}
\caption{Another example: starting with initial values $x_i(0) = i$ for $1 \leq i \leq 100$, the associated system of ODEs approaches the zeros of the Laguerre polynomial $L_{100}$. } 
\end{figure}
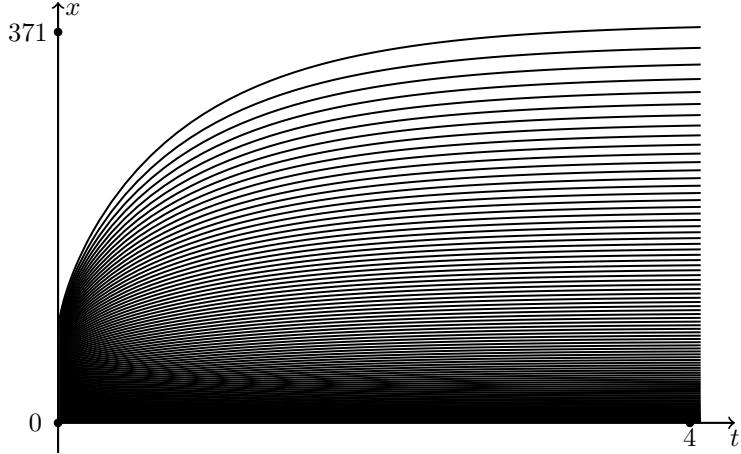
\end{center}
\vspace{-10pt}

\subsection{Proof of Theorem 2.}
\begin{proof} We assume that the solutions $y_0, y_1, \dots, y_{n-1}, y_n$ satisfy $\mbox{deg} (y_i) = i$ for all $0 \leq i \leq n$. We assume that $p, q$ and the domain of definition are such that classical
Sturm-Liouville theory applies. Then the polynomial $y_j$ has exactly $j$ zeros all of which are simple. Suppose  
$$ x_1(0) < x_2(0) < \dots < x_n(0)$$
is given. We define the function
$$ f(x) = \prod_{k=1}^{n}{(x - x_k(0))} \quad \mbox{and write it as} \quad f(x) = \sum_{k=0}^{n}{a_k y_k(x)}$$
for some coefficients $a_0, \dots, a_n$. 
This is possible because of the assumption on the degrees (an upper triangular matrix is invertible).
The strong form of the Sturm Oscillation Theorem, that is not very widely known, states that, as long as not all coefficients $b_i$ vanish, any function of the form
$$ \sum_{k=0}^{n-1}{b_k y_k(x)} \qquad \mbox{has at most}~n-1~\mbox{zeros.}$$
We refer to B\'erard \& Helffer \cite{berard} and L\"utzen \cite{lutzen} for the history of this remarkable Theorem that seems to have been forgotten (\cite{berard} gives rigorous proofs in modern language,  \cite{stein}
gives a quantitative form). Since $f(x)$ has $n$ zeros, the Sturm Oscillation Theorem implies that $a_n \neq 0$. We now define $u(t,x)$ as the solution of the \textit{backward} heat equation 
\begin{align*}
\frac{\partial}{\partial t} u(t,x) &= - \frac{\partial}{\partial x}\left( p(x) \frac{\partial}{\partial x} u(t,x) \right) + q(x) \frac{\partial}{\partial x} u(t,x) \\
u(0,x) &= f(x)
\end{align*}
and observe that, as explained in the proof of Theorem 1, this equation is well-posed for all $t>0$ since it can be rewritten as a linear system of $n+1$ ordinary differential equations. Linearity implies that the solution is given by
$$ u(t,x) = \sum_{k=0}^{n}{a_k e^{\lambda_k t}  y_k(x)}.$$
At the same time, as long as the zeros do not collide, we can write
$$ u(t,x) =  h(t)\prod_{i=1}^{n}{(x - x_i(t))},$$
where $h(t) \neq 0$ for all $t>0$ and the functions $x_i(t)$ satisfy, following the computation done in the proof of Theorem 1 and reversing the sign, the system of ordinary differential equations
$$    \frac{d}{dt} x_i(t) =  -p(x_i) \sum_{k = 1 \atop k \neq i}^{n}{\frac{2}{x_k(t) - x_i(t)}} +  q(x_i(t)) -  p'(x_i(t)) \quad \mbox{for all}~1 \leq i \leq n.$$
It is a property of Sturm-Liouville problems that the number of distinct zeros is nonincreasing in time under the forward heat equation; since we are 
dealing with the backward heat equation, we see that the number of distinct zeros is nondecreasing. Moreover, since we start with a polynomial of degree
$n$ and the solution is always a polynomial of degree at most $n$, this implies that zeros can never collide (and that the solution is always exactly of degree $n$). This shows that
$$ e^{-\lambda_n t} h(t)\prod_{i=1}^{n}{(x - x_i(t))} = a_n y_n(t) +  \sum_{k=0}^{n-1}{a_k e^{(\lambda_k - \lambda_n) t}  y_k(x)}.$$
All zeros of $y_n$ are simple, the Inverse Function Theorem now implies that the zeros of $u(t,x)$ converge to the zeros of $y_n$ exponentially
quickly in $t$. The speed of convergence depends on the size of $\lambda_{n-1} - \lambda_n$ and the size of $y_n'$ at its zeros; the constant in front will depend on the precise values of the coefficients $\left\{a_0, a_1, \dots, a_n\right\}$. If $a_{n-1} = 0$ or more of the leading terms vanish, then convergence would be even faster.
\end{proof}

\end{document}